\documentclass{birkjour}
\usepackage{amsmath}
\usepackage{amssymb}
\usepackage{amsthm}
\usepackage[all]{xy}
\usepackage{enumitem}
\usepackage[T1]{fontenc}
\usepackage{mathtools}

\newtheorem{thm}{Theorem}[section]
\newtheorem{cor}[thm]{Corollary}
\newtheorem{lem}[thm]{Lemma}
\newtheorem{prop}[thm]{Proposition}
\theoremstyle{definition}
\newtheorem{defn}[thm]{Definition}
\theoremstyle{remark}
\newtheorem{rem}[thm]{Remark}

\numberwithin{equation}{section}

\begin{document}

\title{Hamiltonian dynamics and geometry on the two-plectic six-sphere}

\author[Wagner]{Maxime Wagner}
\address{Institut \'Elie Cartan de Lorraine\br
Universit\'e de Lorraine et CNRS\br
57070 Metz, France}
\email{maxime.wagner@univ-lorraine.fr}

\author[Wurzbacher]{Tilmann Wurzbacher}
\address{Institut \'Elie Cartan de Lorraine\br
Universit\'e de Lorraine et CNRS\br
57070 Metz, France}
\email{tilmann.wurzbacher@univ-lorraine.fr}



\begin{abstract}
We study the two-plectic geometry of the six-sphere induced by pulling back a canonical $G_2$-invariant three-form from $\mathbb{R}^7$. Notably we explicitly prove non-flatness of this structure and show that its infinitesimal automorphisms are given by the exceptional Lie algebra $\mathfrak{g}_2$. Several interesting classes of solutions of the dynamical Hamilton-de Donder-Weyl equations with one- and two-dimensional sources are exhibited.
\end{abstract}

\maketitle
\tableofcontents
\clearpage

\section{Introduction}\label{Introduction}

\noindent It is well-known that the six-sphere $\mathbb{S}^6$, viewed as the unit sphere in the imaginary octonions, inherits a non-integrable almost complex structure $J$ from the vector product of octonions. This product, together with the standard scalar product $<,>$ on  $\mathbb{R}^7$, furnishes a closed 3-form $\tilde{\omega}$ as well:
\begin{eqnarray*}
    \tilde{\omega}(u, v, w) := <u \times v, w>.
\end{eqnarray*}

\noindent The inclusion $j: \mathbb{S}^6 \hookrightarrow \mathbb{R}^7$ yields $\omega := j^* \tilde{\omega}$, and both forms are non-degenerate in the sense that contraction with a non-vanishing tangent vector is never zero. (Such 3-forms are called two-plectic.)\\

\noindent Similarly to the case of the two-sphere $\mathbb{S}^2$ that inherits a closed 2-form and an integrable almost complex structure from the quaternions, related in the sense of a K\"ahler structure, the tensors $\omega$ and $J$ enjoy here the following compatibility:
\begin{eqnarray*}
    \omega(Ju, v, w) = \omega(u, Jv, w) = \omega(u, v, Jw) \quad \forall\; u, v, w \in T_p S^6 \text{ and } \forall\; p \in S^6.
\end{eqnarray*}
\noindent Unfortunately, Darboux coordinates are rare on multisymplectic manifolds and the non-integrability of $J$ can be used to show that "flat" coordinates do not exist on $(\mathbb{S}^6,\,\omega)$ either. Nevertheless, the close relation of $\omega$ to the octonions allows to obtain strong constraints on the two-plectic geometry of the six-sphere, notably regarding the automorphism group of $\omega$ that turns out to be the exceptional Lie group $G_2$.\\

\noindent Furthermore, $(\mathbb{S}^6,\,\omega)$ serves here as a case study of Hamiltonian dynamics beyond the symplectic case. Let us recall that, given a closed $(n+1)$-form $\omega$ on a manifold $M$, a couple $(H,\,X)$ in $\Omega^{n{-}k}(M)\oplus \mathfrak{X}^k(M)$ solves the Hamilton-de Donder-Weyl equation in dimension $k$ (or HDW equation) if 
\begin{eqnarray*}
    \iota_X\omega=(-1)^{n+1{-}k}dH.
\end{eqnarray*}
\noindent Similarly, a smooth map $\psi\::\Sigma\rightarrow M$, defined on a $k$-dimensional manifold $\Sigma$ with co-volume $\gamma\in\mathfrak{X}^k(\Sigma)$ solves the dynamical Hamilton-de Donder-Weyl equation in dimension $k$ $\big($with respect to $H\in\Omega^{n{-}k}(M)\big)$ if 
\begin{eqnarray*}
    \forall\; x\in \Sigma,\;\;\iota_{(\psi_*)_x(\gamma_x)}\omega=(-1)^{n+1{-}k}(dH)_{\psi(x)}.
\end{eqnarray*}
\noindent On the two-plectic six-sphere we only have the possibilities $k=1$ or $k=2$, both of which will be studied here.\\

\noindent Let us now proceed to a detailed description of the content of the sections of this article.\\
\\
\noindent In Section 2 we recall the vector product on $\mathbb{R}^7$, viewed as the imaginary octonions, and define the canonical multisymplectic three-form $\widetilde{\omega}$ induced by it. We recall some useful results on $\widetilde{\omega}$ and note that the induced 3-form on $\mathbb{S}^6$ is non-degenerate as well and possesses a $G_2$-invariant potential. Let us underline that in stark contrast to the symplectic case, we exhibit here a compact manifold with a multisymplectic three-form whose cohomology class vanishes. Note that on the other hand the very non-vanishing of the cohomology class of a symplectic form on a compact manifold prevents $\mathbb{S}^6$ from having a symplectic structure.\\

The isomorphy class of a couple $(V,\,\alpha)$ consisting of a finite dimensional vector space $V$ and a non-degenerate $(k+1)$-form $\alpha$ is called the "linear type of $(V,\,\alpha)$". Obviously locally constant linear type is a necessary condition for Darboux coordinates on a multisymplectic manifolds (i.e. local coordinates that render the coefficients of the form constant).\\

\noindent In case $V$ is six-dimensional and $\alpha$ two-plectic there exists only three linear types and we check that $\omega$ on $\mathbb{S}^6$ is everywhere of "complex type", i.e. on each tangent space $T_p\mathbb{S}^6$ there exists a linear almost complex structure $J_p$ and a $\mathbb{C}$-linear 3-form $\Omega_p$ such that $\omega$ is the real part of $\Omega$. For this linear type we know that the above cited algebraic compatibility with the "usual" (non-integrable) $J$ holds. We give a hands-on proof that the existence of local Darboux coordinates for $\omega$ on $\mathbb{S}^6$ implies the existence of a complex analytic atlas inducing $J$, and thus by contradiction of the non-flatness of $\omega$.\\

\noindent In Section 3 we extend the almost complex structure $J$ to an almost Cauchy-Riemann structure on $\mathbb{R}^7\backslash\{0\}$, allowing us to show notably the following identity for vector fields  $X,\,Y,$ and $Z$ on $\mathbb{S}^6$\::
\begin{eqnarray*}
    \omega(X,\,Y,\,Z)=-\frac{1}{4}g\big(N_J(X,\,Y),\,JZ\big),
\end{eqnarray*}
where $N_J$ is the Nijenhuis tensor of $J$ and $g$ the round metric on $\mathbb{S}^6$. We apply this identity to show that the connected component of the automorphism group of $(\mathbb{S}^6,\,\omega)$ is the Lie group $G_2$, sitting inside $SO(7,\,\mathbb{R})$ and thus acting "linearly" on $\mathbb{S}^6$.\\

\noindent We close Section 3 with the surprising observation that the vector space of fundamental vector fields of the $G_2$-action on $\mathbb{S}^6$ is "totally real" in the sense that $\{J(X_{\xi})\;|\;\xi\in\mathfrak{g}_2\}\cap\{X_{\eta}\;|\;\eta\in\mathfrak{g}_2\}=\{0\}$. This result is in stark contrast to the fact that for an integrable almost complex structure on a smooth real manifold the space of of real vector fields preserving this structure is stable under it.\\

\noindent The last section, Section 4, studies Hamiltonian dynamics on the two-plectic six-sphere, more precisely solutions of the dynamical Hamilton-de Donder-Weyl equation with source $\Sigma$ open in $\mathbb{R}^k$ for $k=1$ and $2$, and $\gamma\in \mathfrak{X}^k(\mathbb{R}^k)$ given by the standard co-volume $\frac{\partial}{\partial t^1}\wedge \ldots\wedge \frac{\partial}{\partial t^k}$. The fact that $aut(\mathbb{S}^6,\,\omega)=\mathfrak{g}_2$, shown in Section 3, easily implies that the vector field part of a solution couple of the HDW equation in dimension 1 is a fundamental vector field  of the $G_2$-action on $\mathbb{S}^6$. We give a complete classification of solutions of the HDW equation in this dimension and describe the ensuing "dynamical" solutions $\psi_p^{\xi}\::\mathbb{R}\rightarrow\mathbb{S}^6$ in detail.\\

\noindent In dimension two we exhibit two classes of solutions of the HDW equation. Firstly, a couple $(\xi,\,\eta)$ of commuting elements of $\mathfrak{g}_2$ yields a solution, whose dynamical realisation $\psi_p^{\xi,\eta}$ generalises the above map $\psi_p^{\xi}$ in an obvious manner. We study the image set of these "fields" $\psi_p^{\xi,\eta}$ in some detail.\\
\noindent Finally, we exhibit another surprising difference to Hamiltonian dynamics on a symplectic manifold\:: given a constant function $H$ on a two-plectic manifold $(M,\,\mu)$ there might exist interesting Hamiltonian bi-vector fields $Z$ fulfilling 
\begin{eqnarray*}
    \iota_Z\mu=0=-dH.
\end{eqnarray*}
\noindent In the case of $(\mathbb{S}^6,\,\omega)$ it turns out that for all vector fields $X$ one has 
\begin{eqnarray*}
    \iota_{X\wedge JX}\omega=0,
\end{eqnarray*}
implying that the couple $(0,\,X\wedge JX)$ solves the Hamilton-de Donder-Weyl equation. In case $X$ and $JX$ commute, we obtain a singular $J$-holomorphic foliation whose leaves are the images of solutions of the dynamical version of this equation.\\

\section{Two-plectic structures related to the octonions}\label{Section : Two-plectic structures related to octonions}

We use the octonions to construct canonical $2$-plectic structures on $\mathbb{R}^7$ and $\mathbb{S}^6$ and recall several useful related facts, notably on its automorphisms and on associated almost complex structures. We give a proof of the folkloric result that the $2$-plectic form $\omega$ on $\mathbb{S}^6$ is of "complex" type and show, as the central result of this section that $(\mathbb{S}^6,\,\omega)$ is not flat, i.e. that it does not allow for Darboux coordinates.
\newline
\newline

\noindent Let $\mathbb{O}$ be the $\mathbb{R}$-algebra of octonions defined as follows,
\begin{eqnarray*}
    \mathbb{O}=\Bigg\{x^0e_0+\sum\limits_{i=1}^7x^ie_i,\, x^i\in \mathbb{R} \,\,\Bigg|\,\, e_0=1 \text{ and } e_ie_j=-\delta_i^j+\epsilon_{ijk}e_k \,\,\forall i,j\in \{1,\ldots,7\}\Bigg\}
\end{eqnarray*}
where $\epsilon_{ijk}$ is a completely antisymmetric tensor, equal to $1$ if $ijk= 123$, $145$, $167$, $246$, $275$, $374$, $365$ and equal to $0$ otherwise. We recall the multiplication table of $\mathbb{O}$ (note that $1=e_0$)\::
\begin{center}
    \begin{tabular}{|c||c|c|c|c|c|c|c|}
        \hline
        $\cdot$ & $e_1$ & $e_2$ & $e_3$ & $e_4$ & $e_5$ & $e_6$ & $e_7$ \\
        \hline
        \hline
        $e_1$ & -$1$ & $e_3$ & -$e_2$ & $e_5$ & -$e_4$ & $e_7$ & -$e_6$ \\
        \hline
        $e_2$ & -$e_3$ & -$1$ & $e_1$ & $e_6$ & -$e_7$ & -$e_4$ & $e_5$ \\
        \hline
        $e_3$ & $e_2$ & -$e_1$ & -$1$ & -$e_7$ & -$e_6$ & $e_5$ & $e_4$ \\
        \hline
        $e_4$ & -$e_5$ & -$e_6$ & $e_7$ & -$1$ & $e_1$ & $e_2$ & -$e_3$ \\
        \hline
        $e_5$ & $e_4$ & $e_7$ & $e_6$ & -$e_1$ & -$1$ & -$e_3$ & -$e_7$ \\
        \hline
        $e_6$ & -$e_7$ & $e_4$ & -$e_5$ & -$e_2$ & $e_3$ & -$1$ & $e_1$ \\
        \hline
        $e_7$ & $e_6$ & -$e_5$ & -$e_4$ & $e_3$ & $e_7$ & -$e_1$ & -$1$ \\
        \hline
    \end{tabular}
\end{center}
Either $x^0e_0$ or the coefficient $x^0$ is called the real part and the $x^i$ for $i\in\{1,\ldots,\,7\}$ or $\sum\limits_{i=1}^7x^ie_i$  form the imaginary part of an octonion. As a vector space, we can identify $\mathbb{O}$ with $\mathbb{R}^8$, and the imaginary octonions with $\mathbb{R}^7$. We fix the standard scalar product $<,>$ on $\mathbb{R}^7$ and define a cross product on $\mathbb{R}^7$ via the octonionic multiplication\::
\begin{eqnarray*}
    u\times v:=\frac{1}{2}[u,\,v]=\frac{1}{2}(u\cdot v-v\cdot u).
\end{eqnarray*}

\begin{defn}
    Let $\widetilde{\omega}\in \Omega^3(\mathbb{R}^7)$ be defined by
    \begin{eqnarray*}
        \widetilde{\omega}(u,\,v,\,w)=<u\times v,\,w>
    \end{eqnarray*}
    for $p\in\mathbb{R}^7$ and $u,\,v,\,w\in T_p\mathbb{R}^7$.
\end{defn}

\noindent Direct inspection shows the following result.

\begin{lem}
    The couple $(\mathbb{R}^7,\,\widetilde{\omega})$ is a $2$-plectic manifold.
\end{lem}

\begin{rem}
    Recall that a couple $(M,\,\eta)$ with $\eta$ a closed $(n+1)$-form on $M$ is called an $n$-plectic manifold if for $v$ in $TM$, $\iota_v\eta=0$ implies that $v=0$.
\end{rem}

\begin{defn}
    \raggedright Let $Aut_{\mathbb{R}\text{-}alg}(\mathbb{O})\xhookrightarrow{i} GL(7,\,\mathbb{R})$ be the natural inclusion of the  $\mathbb{R}$-algebra automorphisms of 
    $\mathbb{O}$ given by restriction to the imaginary octonions. Then we define $G_2:=i\big(Aut_{\mathbb{R}\text{-}alg}(\mathbb{O})\big)$.
\end{defn}

\noindent The following results are classical resp. due to R. Bryant (cf. \cite{Bryant1987} and \cite{Bryant2006}).

\begin{lem}\label{Bryant}
Let $\widetilde{\omega}$ and $G_2$ be as above. Then\::

    \begin{enumerate}[label=(\roman*)]
        \item $G_2=\{A\in GL(7,\,\mathbb{R})\,|\, A^*\widetilde{\omega}=\widetilde{\omega}\}$ (the compact real Lie group associated to the exceptional complex Lie algebra 
        $(\mathfrak{g}_2)^{\mathbb{C}}$)\;;
        \item the group $G_2$ is a closed subgroup of $SO(7,\mathbb{R})$\;;
        \item a smooth almost complex structure $J$ is given on $\mathbb{S}^6$ as follows. At each point $p\in\mathbb{S}^6$ one defines \::
        \begin{eqnarray*}
    J_p\:: T_p\mathbb{S}^6  \rightarrow  T_p\mathbb{S}^6,\hspace{0.3cm} u  \mapsto  J_p(u):=p\times u \, ,
\end{eqnarray*}
where $p\times u$ is the cross product between $p$ and $u$\;;
        \item $G_2$ acts transitively on $\mathbb{S}^6$, preserves $J$ and the isotropy subgroups are isomorphic to $SU(3)$\;;
        \item defining $j\:: \mathbb{S}^6\hookrightarrow \mathbb{R}^7$ as the canonical inclusion and $\omega:=j^*\widetilde{\omega}$, $G_2$ preserves $\omega$\;;
        \item the couple $(\mathbb{S}^6,\,\omega)$ is a $2$-plectic manifold\;;
        \item the form $\omega$ has a $G_2$-invariant potential $\theta$, e.g. $j^*\big(\frac{1}{3}\iota_E\widetilde{\omega}\big)$, where $E$ is the Euler vector field on $\mathbb{R}^7$.
    \end{enumerate}
\end{lem}

\begin{rem}\label{strong non-integrability}
    It is well known that the above defined almost complex structure $J$ is not integrable. The first proofs of this fact were given by B. Eckmann and A. Fr\"olicher (\cite{Eckmann&Frolicher}) resp. C. Ehresmann and P. Liebermann (\cite{Ehresmann&Libermann}). A stronger result  can be found in an article of B. Kruglikov  (\cite{Kruglikov}) showing the non-degeneration of the Nijenhuis tensor. More precisely it is there shown that (using $J$ as multiplication by $\sqrt{-1}$) 
    \begin{eqnarray*}
        N_J\:: \Lambda^2T\mathbb{S}^6\rightarrow T\mathbb{S}^6
    \end{eqnarray*}
    is a $\mathbb{C}$-linear vector bundle isomorphism.
\end{rem}

\begin{defn}
    A $k$-plectic \textit{linear type} is an isomorphism class of $k$-plectic vector spaces (where $(V,\,\omega)$ is \textit{isomorphic} to $(V',\,\omega')$ if and only if there exists a linear isomorphism $T\::V\rightarrow V'$ such that $T^*\omega'=\omega$).
\end{defn}

\noindent \raggedright To determine the linear type of $\omega$ on $\mathbb{S}^6$, we use tools developped by M. Pan\'{a}k and J. Van\v{z}ura (\cite{Panak&Vanzura}).

\begin{defn}
    For $V$ a real vector space of dimension $6$ and $\alpha\in\Lambda^3V^*$, $\Delta(\alpha)$ is defined by $\{v\in V\;|\; \iota_v\alpha\wedge\iota_v\alpha=0\}$.
\end{defn}

\noindent Using this definition, we can reformulate a result of B. Capdevielle.

\begin{prop}[\cite{Capdevielle}]
    Let $V$ be a real vector space of dimension $6$. Then there exist three orbits for the action of $GL_{\mathbb{R}}(V)$ on non-degenerate forms $\alpha$ in $\Lambda^3V^*$ determined by $\Delta$ as follow\::
    \begin{enumerate}[label=(\roman*)]
        \item if $\Delta(\alpha)=V^a\cup V^b$ with $V^a$ and $V^b$ of dimension three and $V^a\cap V^b=\{0\}$, $\alpha$ is of product type\;;
        \item if $\Delta(\alpha)=\{0\}$, $\alpha$ is of complex type\;;
        \item if $\Delta(\alpha)$ is a three-dimensional subspace, $\alpha$ is of tangent type.
    \end{enumerate}
\end{prop}

\begin{rem}
    We use this proposition as the definition of the three types, but remark that $\alpha$ is of complex type if and only if $V$ has a complex structure and $\Omega\in \underset{\mathbb{C}}{\Lambda}^3V^*$ such that $\omega=\Re(\Omega)$.
\end{rem}

\begin{cor}\label{corollaire GL-orbites}
    Let $V$ be a real vector space of dimension $6$, then $GL_{\mathbb{R}}(V)$ acts transitively on the set of non-degenerate forms $\alpha\in\Lambda^3V^*$ of complex type.
\end{cor}

\begin{lem}
    Let $p\in \mathbb{S}^6$ and $\omega$ as in Lemma \ref{Bryant}, then we have 
    $\Delta(\omega_p)=\{v\in T_p\mathbb{S}^6\;|\; \iota_v\omega_p\wedge\iota_v\omega_p=0\}=\{0\}$, i.e. $\omega_p$ is of complex type.
\end{lem}

\begin{proof}
    By the transitivity of the $G_2$-action on $\mathbb{S}^6$, we can assume $p=N=(1,\,0,\ldots,0)$. Then $T_N\mathbb{S}^6=e_1^{\perp}\subset \mathbb{R}^7$ and \begin{eqnarray*}
        \omega_N=dx^2\wedge(dx^4\wedge dx^6-dx^5\wedge dx^7)-dx^3\wedge (dx^4\wedge dx^7+dx^5\wedge dx^6). 
    \end{eqnarray*}
    With $v=\sum\limits_{i=1}^6 v^i\partial_{x^{i+1}}$ denoting a tangent vector in $T_N\mathbb{S}^6$, one has
    \begin{eqnarray*}
        \iota_v\omega_N=v^1(dx^4\wedge dx^6-dx^5\wedge dx^7)+v^2(-dx^4\wedge dx^7-dx^5\wedge dx^6)\\
        +v^3(-dx^2\wedge dx^6+dx^3\wedge dx^7)+v^4(dx^2\wedge dx^7+dx^3\wedge dx^6)\\
        +v^5(dx^2\wedge dx^4-dx^3\wedge dx^5)+v^6(-dx^2\wedge dx^5-dx^3\wedge dx^4).
    \end{eqnarray*}
    Thus we obtain $(\iota_v\omega_N)\wedge \omega_N=$
    \begin{eqnarray*}
        & & v^1(-dx^4\wedge dx^6\wedge dx^2\wedge dx^5\wedge dx^7-dx^5\wedge dx^7\wedge dx^2\wedge dx^4\wedge dx^6)\\
        & & +v^2(dx^4\wedge dx^7\wedge dx^3\wedge dx^5\wedge dx^6+dx^5\wedge dx^6\wedge dx^3\wedge dx^4\wedge dx^7)\\
        & & +v^3(dx^2\wedge dx^6\wedge dx^3\wedge dx^4\wedge dx^7+dx^3\wedge dx^7\wedge dx^2\wedge dx^4\wedge dx^6)\\
        & & +v^4(-dx^2\wedge dx^7\wedge dx^3\wedge dx^5\wedge dx^6-dx^3\wedge dx^6\wedge dx^2\wedge dx^5\wedge dx^7)\\
        & & +v^5(dx^2\wedge dx^4\wedge dx^3\wedge dx^5\wedge dx^6-dx^3\wedge dx^5\wedge dx^2\wedge dx^4\wedge dx^6)\\
        & & +v^6(dx^2\wedge dx^5\wedge dx^3\wedge dx^4\wedge dx^7+dx^3\wedge dx^4\wedge dx^2\wedge dx^5\wedge dx^7)\\
        & = & 2 v^1 dx^2\wedge dx^4\wedge dx^5\wedge dx^6\wedge dx^7 +2 v^2 dx^3\wedge dx^4\wedge dx^5\wedge dx^6\wedge dx^7 \\
        & & +2 v^3 dx^2\wedge dx^3\wedge dx^4\wedge dx^6\wedge dx^7 +2 v^4 dx^2\wedge dx^3\wedge dx^5\wedge dx^6\wedge dx^7 \\
        & & +2 v^5 dx^2\wedge dx^3\wedge dx^4\wedge dx^5\wedge dx^6 +2 v^6 dx^2\wedge dx^3\wedge dx^4\wedge dx^5\wedge dx^7.
    \end{eqnarray*}
    Let $\nu=dx^2\wedge dx^3\wedge dx^4\wedge dx^5\wedge dx^6\wedge dx^7$ be the canonical volume form on  $T_N\mathbb{S}^6\cong \{x^1\!=\!0\}\subset \mathbb{R}^7$. Then there exists a unique tangent vector $w_v$ such that $(\iota_v\omega_N)\wedge \omega_N=\iota_{w_v}\nu$. This tangent vector is $w_v=2(v^2\partial_{x^2}-v^1\partial_{x^3}+v^4\partial_{x^4}-v^3\partial_{x^5}+v^6\partial_{x^6}-v^5\partial_{x^7})$. Moreover
    \begin{eqnarray*}
        \iota_v\big(\iota_v\omega_N\wedge \omega_N\big) & =  &\iota_v\iota_v\omega_N\wedge \omega_N+\iota_v\omega_N\wedge \iota_v\omega_N=\iota_v\omega_N\wedge \iota_v\omega_N\\
        & = & \iota_v\iota_{w_v}\nu.
    \end{eqnarray*}
    So, $\iota_v\omega_N\wedge \iota_v\omega_N=0$ is equivalent to $v\wedge w_v=0$. But explicit calculations show that
    \begin{eqnarray*}
        0=v\wedge w_v \text{ implies}  
 & -2\big((v^1)^2+(v^2)^2\big)\partial_{x^2}\wedge \partial_{x^3} -2 \big((v^3)^2+(v^4)^2\big)\partial_{x^4}\wedge \partial_{x^5}\\
 & -2 \big((v^5)^2+(v^6)^2\big)\partial_{x^6}\wedge \partial_{x^7}=0.
    \end{eqnarray*}
    Thus $(v^1)^2+(v^2)^2=0\;, (v^3)^2+(v^4)^2 =0$
    and $(v^5)^2+(v^6)^2=0$, and therefore $v^1=v^2=v^3=v^4=v^5=v^6=0$.\\
    Thus we have $\Delta(\omega_N)=\{0\}$.
\end{proof}

\begin{prop}[\cite{Panak&Vanzura}]\label{Lemme technique}
    Let $V$ be a real vector space of dimension $6$ and $\alpha\in\Lambda^3V^*$ non-degenerate with $\Delta(\alpha)=\{0\}$. Then 
    \begin{enumerate}[label=(\roman*)]
        \item there exists, up to a sign, a unique almost complex structure $J$ on $V$ satisfying
        $$\alpha(Ju,\,v,\,w)=\alpha(u,\,Jv,\,w)=\alpha(u,\,v,\,Jw) \,\,\, \forall u,\,v,\,w\in V   \;; $$
        \item$\ker(\iota_v\alpha)=\big(\big(v,\,J(v)\big)\big)_{\mathbb{R}}$  for all $v\in V\backslash\{0\}$.
    \end{enumerate}
\end{prop}

\begin{defn}
    Let $(M,\,\omega)$ be a $k$-plectic manifold and $p$ in $M$. We say $M$ has \textit{Darboux coordinates} near $p$ if it exists a neighbourhood $U$ of $p$ with coordinates $u^1,\,\ldots,\,u^n$ and $a_{i_1,\,\ldots,\,i_{k+1}}\in \mathbb{R}$ such that
    \begin{eqnarray*}
        \omega_{|U}=\sum_{1\leq i_1<\,\ldots <i_{k+1}\leq n}a_{i_1,\,\ldots,\,i_{k+1}}du^{i_1}\wedge \ldots \wedge du^{i_{k+1}}.
    \end{eqnarray*}
\end{defn}

\noindent \raggedright The next result is similar to Theorem $4.11$ in \cite{Ryvkin&Wurzbacher} and to Corollary $15$ in \cite{Panak&Vanzura}. In contrast to these two articles, we explicitely show that the local existence of Darboux coordinates for $\omega$, e.g. near the north pole, implies the existence of a complex analytic atlas on $(\mathbb{S}^6,\,J)$.

\begin{thm}
    On $(\mathbb{S}^6,\,\omega)$ there exist no Darboux coordinates.
\end{thm}

\begin{proof}
    Suppose that near the point $N$, there exists a Darboux chart $(U,\,\widetilde{\varphi})$, i.e. a form on $\mathbb{R}^6$,
    \begin{eqnarray*}
        \widehat{\omega}=\sum_{1\leq i_1<i_2 <i_3\leq 6}a_{i_1,\,i_2,\,i_3}du^{i_1}\wedge du^{i_2}\wedge du^{i_3}
    \end{eqnarray*}
    with $a_{i_1,\,i_2,\,i_3}\in\mathbb{R}$ and a chart $(U,\,\widetilde{\varphi})$ such that $\widetilde{\varphi}^*\widehat{\omega}=\omega_{|U}$.
    \newline
    
    By Corollary \ref{corollaire GL-orbites}, $GL(6,\,\mathbb{R})$ acts transitively on $2$-plectic forms $\alpha$ in $\Lambda^3(\mathbb{R}^6)^*$ with $\Delta(\alpha)=\{0\}$. 
    Since 
    $\check{\omega}=dx^1\wedge dx^2\wedge dx^3-dx^1\wedge dy^2\wedge dy^3-dy^1\wedge dx^2\wedge dy^3-dy^1\wedge dy^2\wedge dx^3$ is $2$-plectic with
    $\Delta(\check{\omega})=\{0\}$,
    there exists $A\in GL(6,\,\mathbb{R})$ such that $T_A^*\check{\omega}=\hat{\omega}$.
    \newline
    
    We denote $T_A\circ \widetilde{\varphi}$ by $\varphi$. Observe that on $\varphi(U)\subset \mathbb{R}^6$ we have complex analytic coordinates $z^k=x^k+\sqrt{-1}y^k$ $(k=1,\,2,\,3)$ such that $\check{\omega}=\Re(dz^1\wedge dz^2\wedge dz^3)$ and the induced almost complex structure fulfills $J^{\mathbb{R}^6}\Big(\frac{\partial}{\partial x^k}\Big)=\frac{\partial}{\partial y^k}$ for $k=1,\,2,\,3$. By Lemma \ref{Lemme technique} there exists an almost complex structure $J$, unique up to a sign, such that for all $p\in\mathbb{S}^6$ and $u,\,v,\,w\in T_p\mathbb{S}^6$, one has the following compatibility $\omega_p(J_pu,\,v,\,w)=\omega_p(u,\,J_pv,\,w)=\omega_p(u,\,v,\,J_pw)$. Since $(\varphi^{-1})^*\omega=\check{\omega}$, and $\varphi_*\circ J\circ (\varphi^{-1})_*$ as well as $J^{\mathbb{R}^6}$ enjoy the same compatibility with $\check{\omega}$, we have $J^{\mathbb{R}^6}=\pm \varphi_*\circ J\circ (\varphi^{-1})_*$. We conclude that, upon possibly replacing $z^j$ by $\overline{z^j}$, $\varphi_*\circ J \circ (\varphi^{-1})_*=J^{\mathbb{R}^6}$ and $(\varphi^{-1})^*\omega=\Re(dz^1\wedge dz^2\wedge dz^3)$.
    \newline
    
    We define $U_g=\vartheta_g(U)$ for $g\in G_2$ and $\vartheta_g : \mathbb{S}^6 \rightarrow \mathbb{S}^6$, the transitive action of $G_2$ on $\mathbb{S}^6$. Then we have 
    $$\vartheta^*_g(\omega_{|\vartheta_g(U)})=\omega_{|U} \text{  and  } (\vartheta_g)_*(J_{|U})=J_{|\vartheta_g(U)},$$
    and we set $\varphi_{_N}:=\varphi$,  respectively $\varphi_{\vartheta_g(N)}:=\varphi_{_N}\circ\vartheta_g^{-1}$, defined by $\varphi_{_N}$ near $N$, respectively $\vartheta_g(N)$. Noting  $\Psi=\varphi_{_N}\circ\varphi_{\vartheta_g(N)}^{-1}$, we have the following equality for $p\in \varphi_{_N}(U)\cap\varphi_{\vartheta_g(N)}(U_g)$\:: 
    \begin{eqnarray*}
        & & T_{p}\Psi\circ J_p^{\mathbb{R}^6}=  T_p\big(\varphi_{_N}\circ\varphi_{\vartheta_g(N)}^{-1}\big)\circ J_p^{\mathbb{R}^6}\\
        & = & T_p\big(\varphi_{_N}\circ\varphi_{\vartheta_g(N)}^{-1}\big)\circ\big(T_{\varphi_{\vartheta_g(N)}^{-1}(p)}\varphi_{\vartheta_g(N)}\circ J_{\varphi^{-1}_{\vartheta_g(N)}(p)}\circ T_p\varphi_{\vartheta_g(N)}^{-1}\big)\\
        & = & \big(T_{\varphi_{\vartheta_g(N)}^{-1}(p)}\varphi_{_N}\circ T_p\varphi_{\vartheta_g(N)}^{-1}\big)\circ T_{\varphi_{\vartheta_g(N)}^{-1}(p)}\varphi_{\vartheta_g(N)}\circ J_{\varphi^{-1}_{\vartheta_g(N)}(p)}\circ T_p\varphi_{\vartheta_g(N)}^{-1}\\
        & = & T_{\varphi_{\vartheta_g(N)}^{-1}(p)}\varphi_{_N}\circ\big(T_p\varphi_{\vartheta_g(N)}^{-1}\circ T_{\varphi_{\vartheta_g(N)}^{-1}(p)}\varphi_{\vartheta_g(N)}\big)\circ J_{\varphi^{-1}_{\vartheta_g(N)}(p)}\circ T_p\varphi_{\vartheta_g(N)}^{-1}\\ 
    \end{eqnarray*}
    \begin{eqnarray*}
        & = & T_{\varphi_{\vartheta_g(N)}^{-1}(p)}\varphi_{_N}\circ J_{\varphi^{-1}_{\vartheta_g(N)}(p)}\circ T_p\varphi_{\vartheta_g(N)}^{-1}\\
        & = & T_{\varphi_{\vartheta_g(N)}^{-1}(p)}\varphi_{_N}\circ J_{\varphi^{-1}_{\vartheta_g(N)}(p)}\circ \big(T_{\varphi_{_N}\circ\varphi_{\vartheta_g(N)^{-1}(p)}}\varphi_{_N}^{-1}\circ T_{\varphi_{\vartheta_g(N)}^{-1}(p)}\varphi_{_N}\big)\circ T_p\varphi_{\vartheta_g(N)}^{-1}\\
        & = & \big(T_{\varphi_{\vartheta_g(N)}^{-1}(p)}\varphi_{_N}\circ J_{\varphi^{-1}_{\vartheta_g(N)}(p)}\circ T_{\varphi_{_N}\circ\varphi_{\vartheta_g(N)^{-1}(p)}}\varphi_{_N}^{-1}\big)\circ \big(T_{\varphi_{\vartheta_g(N)}^{-1}(p)}\varphi_{_N}\circ T_p\varphi_{\vartheta_g(N)}^{-1}\big)\\
        & = & J^{\mathbb{R}^6}_{\varphi_{_N}\circ\varphi_{\vartheta_g(N)}^{-1}(p)}\circ T_p\big(\varphi_{_N}\circ\varphi_{\vartheta_g(N)}^{-1}\big)\\
        & = & J^{\mathbb{R}^6}_{\Psi(p)}\circ T_p\Psi.
    \end{eqnarray*}
    Thus $\Psi$ is a complex analytic change of coordinates and one obtains a complex analytic atlas inducing $J$. But since $J$ is not integrable, we are led to a contradiction, implying that there are no Darboux coordinates near $N$, and thus, by the transitivity of the $G_2$-action, nowhere on $\mathbb{S}^6$.
\end{proof}

\section{Automorphisms of the two-plectic six-sphere}\label{Section : Automorphisms of the 2-plectic six-sphere}

In this section we show that the connected component of the group of diffeomorphisms of $\mathbb{S}^6$ preserving the $2$-plectic form $\omega$ is finite dimensional, more precisely equal to the (compact) exceptional Lie group $G_2$. We obtain this result by a formula relating $\omega$ to the round metric and the "usual" almost complex structure $J$ on $\mathbb{S}^6$ and a sequence of six lemmata related to $J$ and its extension to an almost CR stucture on $\mathbb{R}^7\backslash\{0\}$.

\begin{lem}
    Let $X\in\mathfrak{X}^1(\mathbb{S}^6)$ be a vector field, then 
    \begin{eqnarray*}
        (\mathcal{L}_XJ)\circ J=-J\circ(\mathcal{L}_XJ).
    \end{eqnarray*}
\end{lem}

\begin{proof}
    \begin{eqnarray*}
        0= -\mathcal{L}_X(\text{id}_{T\mathbb{S}^6}) = \mathcal{L}_X(J\circ J)=(\mathcal{L}_X J)\circ J+ J\circ (\mathcal{L}_X J).
    \end{eqnarray*}
\end{proof}

\begin{lem}
    Let $X,\,U,\,V,\,W\in\mathfrak{X}^1(\mathbb{S}^6)$ be vector fields such that $\mathcal{L}_X\omega=0$, then we have
    \begin{eqnarray*}
        \omega\big((\mathcal{L}_XJ)U,\,V,\,W\big)=\omega\big(U,\,(\mathcal{L}_XJ)V,\,W\big)=\omega\big(U,\,V,\,(\mathcal{L}_XJ)W\big).
    \end{eqnarray*}
\end{lem}

\begin{proof}
    We have
    \begin{eqnarray*}
        (\mathcal{L}_X\omega)(JU,\,V,\,W) & \!\!\!\!=\!\!\!\! & X\big(\omega(JU,\,V,\,W)\big)+\omega\big((\mathcal{L}_XJ)U,\,V,\,W\big)+\omega\big(J\mathcal{L}_XU,\,V,\,W\big)\\
        & & +\omega\big(JU,\,\mathcal{L}_XV,\,W\big)+\omega\big(JU,\,V,\,\mathcal{L}_XW\big)
    \end{eqnarray*}
    and 
    \begin{eqnarray*}
        (\mathcal{L}_X\omega)(U,\,JV,\,W) & \!\!\!\!=\!\!\!\! & X\big(\omega(U,\,JV,\,W)\big)+\omega\big(\mathcal{L}_XU,\,JV,\,W\big)+\omega\big(U,\,J\mathcal{L}_XV,\,W\big)\\
        & & +\omega\big(U,\,(\mathcal{L}_XJ)V,\,W\big)+\omega\big(U,\,JV,\,\mathcal{L}_XW\big).
    \end{eqnarray*}
    Since $\mathcal{L}_X\omega=0$, Proposition \ref{Lemme technique} implies 
    \begin{eqnarray*}
        0 & = & (\mathcal{L}_X\omega)(JU,\,V,\,W)-(\mathcal{L}_X\omega)(U,\,JV,\,W)\\
        & = & \omega\big((\mathcal{L}_XJ)U,\,V,\,W\big)-\omega\big(U,\,(\mathcal{L}_XJ)V,\,W\big).
    \end{eqnarray*}
    We conclude that $\omega\big((\mathcal{L}_XJ)U,\,V,\,W\big)=\omega\big(U,\,(\mathcal{L}_XJ)V,\,W\big)$. The second equality follows from the skew-symmetry of $\omega$.
\end{proof}

\begin{lem}\label{omega et J}
    Let $X \in \mathfrak{X}^1(\mathbb{S}^6)$ be a vector field, then 
    \begin{eqnarray*}
        \mathcal{L}_X\omega=0\text{ implies } \mathcal{L}_XJ=0.
    \end{eqnarray*}
\end{lem}

\begin{proof}
    With the previous lemmata one obtains
    \begin{eqnarray*}
        -\omega\big(J(\mathcal{L}_XJ)U,\,V,\,W\big) & = & \omega\big((\mathcal{L}_XJ)JU,\,V,\,W\big)\\
        & = & \omega\big(JU,(\mathcal{L}_XJ)V,\,W\big)\\
        & = & \omega\big(U,(\mathcal{L}_XJ)V,\,JW\big)\\
        & = & \omega\big((\mathcal{L}_XJ)U,\,V,\,JW\big)\\
        & = & \omega\big(J(\mathcal{L}_XJ)U,\,V,\,W\big).
    \end{eqnarray*}
    We conclude that $\mathcal{L}_XJ=0$.
\end{proof}

\begin{lem}\label{NJ G2-invariance}
    Let $X\in\mathfrak{X}^1(\mathbb{S}^6)$ be a vector field, then
    \begin{eqnarray*}
        \mathcal{L}_XJ=0\text{ implies } \mathcal{L}_X(N_J)=0.
    \end{eqnarray*}
\end{lem}

\begin{proof}
    First observe that if $\mathcal{L}_XJ=0$, then for all $Y\in\mathfrak{X}^1(\mathbb{S}^6)$
    \begin{eqnarray*}
        \mathcal{L}_X(JY)=J\mathcal{L}_XY
    \end{eqnarray*}
    \begin{eqnarray*}
        \text{ and }\,\, \mathcal{L}_X(N_J)=0 \,\,\text{ if and only if }\,\, \mathcal{L}_X\big(N_J(Y,\,Z)\big)=N_J\big(\mathcal{L}_X(Y),\,Z\big)+N_J\big(Y,\,\mathcal{L}_X(Z)\big).
    \end{eqnarray*}
    We compute $\mathcal{L}_X\big(N_J(Y,\,Z)\big) =$
    \begin{eqnarray*}
        & & \mathcal{L}_X\big([JY,\,JZ]-J[JY,\,Z]-J[Y,\,JZ]-[Y,\,Z]\big)\\
        & = & \mathcal{L}_X\big([JY,\,JZ]\big)-\mathcal{L}_X\big(J[JY,\,Z]\big)-\mathcal{L}_X\big(J[Y,\,JZ]\big)-\mathcal{L}_X\big([Y,\,Z]\big)\\
        & = & [X,\,[JY,\,JZ]]-J(\mathcal{L}_X[JY,\,Z])-J(\mathcal{L}_X[Y,\,JZ])-[X,\,[Y,\,Z]]\\
        & = & [X,\,[JY,\,JZ]]-J([X,\,[JY,\,Z]])-J([X,\,[Y,\,JZ]])-[X,\,[Y,\,Z]]\\
        & = & [[X,\,JY],\,JZ]]-J([[X,\,JY],\,Z])-J([[X,\,Y],\,JZ])-[[X,\,Y],\,Z]\\
        & & +[JY,\,[X,\,JZ]]-J([JY,\,[X,\,Z]])-J([Y,\,[X,\,JZ]])-[Y,\,[X,\,Z]]\\
        & = & [\mathcal{L}_X(JY),\,JZ]-J([\mathcal{L}_X(JY),\,Z])-J([\mathcal{L}_XY,\,JZ])-[\mathcal{L}_XY,\,Z]\\
        & & +[JY,\,\mathcal{L}_X(JZ)]-J([JY,\,\mathcal{L}_XZ])-J([Y,\,\mathcal{L}_X(JZ)])-[Y,\,\mathcal{L}_XZ]\\
        & = & [J\mathcal{L}_XY,\,JZ]-J([J\mathcal{L}_XY,\,Z])-J([\mathcal{L}_XY,\,JZ])-[\mathcal{L}_XY,\,Z]+\\
        & & [JY,\,J\mathcal{L}_XZ]-J([JY,\,\mathcal{L}_XZ])-J([Y,\,J\mathcal{L}_XZ])-[Y,\,\mathcal{L}_XZ]\\
        & = & N_J\big(\mathcal{L}_X(Y),\,Z\big)+N_J\big(Y,\,\mathcal{L}_X(Z)\big)
    \end{eqnarray*}
    and thus $\mathcal{L}_X(N_J)=0$.
\end{proof}

\begin{defn}
    Let $\tilde{p}=(r,\,\theta)\in \mathbb{R}^*_+\times\mathbb{S}^6=\mathbb{R}^7\backslash\{0\}$ a non-zero point of $\mathbb{R}^7$. Observe that $T_{\tilde{p}}(\mathbb{R}^7\backslash\{0\})=\tilde{p}^{\perp}\oplus \mathbb{R}\tilde{p}$. For a vector $\tilde{v}$ of the previous tangent space, we note the respective components $\tilde{v}^{\perp}$ et $\tilde{v}^{\parallel}$. Define
    \begin{eqnarray*}
        \widetilde{J}_{\tilde{p}}: T_{\tilde{p}}\mathbb{R}^7\backslash\{0\} & \rightarrow  & T_{\tilde{p}}\mathbb{R}^7\backslash\{0\}\\
        \tilde{v} & \mapsto & \frac{\tilde{p}}{||\widetilde{p}||}\times\tilde{v}.
    \end{eqnarray*}
\end{defn}

\begin{rem}
    Observe that $\widetilde{J}_{\tilde{p}}|_{\mathbb{R}\tilde{p}}=0$, $\widetilde{J}_{\tilde{p}}$ preserves $\tilde{p}^{\perp}$ and fulfills $(\widetilde{J}_p)^2|_{\tilde{p}^{\perp}}=-\text{id}_{\tilde{p}^{\perp}}$. Then $\widetilde{J}$ is an almost CR structure on $\mathbb{R}^7\backslash\{0\}$ (compare, e.g. \cite{Mizner} for this notion).
\end{rem}

\begin{lem}
    Let $S\subset M$ be an embedded submanifold, $\widetilde{L}$ a $(1,1)$-tensor on $M$ such that $\widetilde{L}(TS)\subset TS$. Then the Nijenhuis tensor $N_L$ of $L:=\widetilde{L}|_{TS}$ equals the Nijenhuis tensor $N_{\widetilde{L}}$ of $\widetilde{L}$ restricted to $TS$.
\end{lem}

\begin{proof}
    The result follows easily from the following fact : let $\widetilde{X},\,\widetilde{Y}\in\mathfrak{X}(M)$ be vector fields that are tangent to $S$. Then, for all $p$ in $S,\, [\widetilde{X},\,\widetilde{Y}](p)=[X,\,Y](p)$, where $X=\widetilde{X}|_S$, $Y=\widetilde{Y}|_S$.
\end{proof}

\begin{lem}\label{J et NJ}
    Let $\widetilde{J}$ be the above almost CR structure on $\mathbb{R}^7\backslash\{0\}$. Then 
    \begin{enumerate}[label=(\roman*)]
        \item for all $p\in\mathbb{S}^6$, $\widetilde{J}_p=J_p$ on $T_p\mathbb{S}^6=p^{\perp}$\;;
        \item for all $p\in\mathbb{S}^6$, $(N_{\widetilde{J}})_p=(N_J)_p$ on $T_p\mathbb{S}^6\times T_p\mathbb{S}^6$\;;
        \item for all $p\in\mathbb{S}^6$ and $u_p,\,v_p\in T_p\mathbb{S}^6\subset T_p(\mathbb{R}^7\backslash\{0\})=\mathbb{R}^7$, $N_J(u_p,\,v_p)=-4\widetilde{J}_p(u_p\times v_p)$.
    \end{enumerate}
\end{lem}

\begin{proof}
    \begin{enumerate}[label=(\roman*)]
        \item is obvious.
        \item follows directly from the preceding lemma.
        \item Let $\widetilde{p}=p=N=(1,\,0,\,\ldots,\,0)\in\mathbb{S}^6\subset \mathbb{R}^7\backslash\{0\}$ and note that all elements of $T_p\mathbb{S}^6\subset T_p\mathbb{R}^7\backslash\{0\}$ are $\mathbb{R}$-linear combinations of coordinates vector fields $\frac{\partial}{\partial x^k}=\partial_k$ with $k\geq 2$ on $\mathbb{R}^7\backslash\{0\}$. We obtain the result for such vector fields by direct verification, e.g. for $k=3$, $l=5$ one has
        \begin{eqnarray*}
             \big(N_{\widetilde{J}}(\partial_{x^3},\,\partial_{x^5})\big)_p = \Big([\widetilde{J}(\partial_{x^3}),\,\widetilde{J}(\partial_{x^5})]-\widetilde{J}[\widetilde{J}(\partial_{x^3}),\,\partial_{x^5}]-\widetilde{J}[\partial_{x^3},\,\widetilde{J}(\partial_{x^5})]-[\partial_{x^3},\,\partial_{x^5}]\Big)_p.
        \end{eqnarray*}
        Since 
        \begin{eqnarray*}
            \Big([\widetilde{J}(\partial_{x^3}),\,\widetilde{J}(\partial_{x^5})]\Big)_p & = & \Big(\big[\frac{1}{r}(x\times \partial_{x^3}),\frac{1}{r}(\,x \times \partial_{x^5})\big]\Big)_p\\
            & = & \Big(\frac{1}{r^2}\big[x\times\partial_{x^3},\,x\times\partial_{x^5}\big]+\frac{1}{r}(x\times\partial_{x^3})(r^{-1})(x\times\partial_{x^5})\\
            & & -\frac{1}{r}(x\times\partial_{x^5})(r^{-1})(x\times\partial_{x^3})\Big)_p\\
            & = & \Big(\big[x\times\partial_{x^3},\,x\times\partial_{x^5}\big]\Big)_p\\
            & = & \Big(\sum\limits_{i,\,j=1}^7x^i(\partial_{x^i}\times\partial_{x^3})(x^j)(\partial_{x^j}\times\partial_{x^5})\\
            & & -x^j(\partial_{x^j}\times\partial_{x^5})(x^i)(\partial_{x^i}\times\partial_{x^3})\Big)_p\\
            & = & (-\partial_{x^2}\times\partial_{x^5}+\partial_{x^4}\times\partial_{x^3})_p\\
            & = & (2\partial_{x^7})_p,
        \end{eqnarray*}
        and  
        \begin{eqnarray*}
            \Big([\widetilde{J}(\partial_{x^3}),\,\partial_{x^5}]\Big)_p & = & \Big(\big[\frac{1}{r}(x\times\partial_{x^3}),\,\partial_{x^5}\big]\Big)_p\\
            & = & \Big(-\partial_{x^5}(r^{-1})(x\times\partial_{x^3})+\frac{1}{r}[x\times\partial_{x^3},\,\partial_{x^5}]\Big)_p\\
            & = & \Big(\frac{1}{r}\sum\limits_{i=1}^7[x^i(\partial_{x^i}\times\partial_{x^3}),\,\partial_{x^5}]\Big)_p\\
            & = & \Big([\partial_{x^1}\times\partial_{x^3},\,\partial_{x^5}]-\partial_{x^5}\times\partial_{x^3}\Big)_p\\
            & = & -(\partial_{x^6})_p,
        \end{eqnarray*}
        $\Big(\widetilde{J}[\widetilde{J}(\partial_{x^3}),\,\partial_{x^5}]\Big)_p=-(\partial_{x^7})_p$. A similar calculation yields $\Big(\widetilde{J}[\partial_{x^3},\,\widetilde{J}(\partial_{x^5})]\Big)_p=-(\partial_{x^7})_p$. Then we have
        \begin{eqnarray*}
            \big(N_J(\partial_{x^3},\,\partial_{x^5})\big)_p=4(\partial_{x^7})_p=-4\widetilde{J}_p(\partial_{x^3}\times\partial_{x^5}).
        \end{eqnarray*}
        We can verify this equality on all vector field generated by elements of the canonical base at the point $N$. Then, by the transitivity of the action of $G_2$ on the six-sphere and the invariance of $N_J$ under this action by Lemma \ref{NJ G2-invariance}, we can extend this result to the entire six-sphere.
    \end{enumerate}
\end{proof}

\noindent The next proposition is motivated by formulas obtained in the articles \cite{Agricola&Borowka&Friedrich} and \cite{Verbitsky}.

\begin{prop}\label{omega=g(NJ...)}
    Let $g$ be the round Riemaniann metric on the six-sphere and $N_J$ the Nijenhuis tensor associated to the almost complex structure $J$ coming from the octonions. Then, for $X,\,Y,\,Z\in\mathfrak{X}^1(\mathbb{S}^6)$, we have 
    \begin{eqnarray*}
        \omega(X,\,Y,\,Z)=-\frac{1}{4}g\big(N_J(X,\,Y)\,,JZ\big).
    \end{eqnarray*}.
\end{prop}

\begin{proof}
    For $\tilde{p}\in \mathbb{R}^7\backslash\{0\}$ and $\tilde{v}\in T_{\tilde{p}}\mathbb{R}^7\backslash\{0\}$, we will use the notation $\tilde{v}^{\perp}$ et $\tilde{v}^{\parallel}$ defined before. Moreover, we will use the notation $\langle\, , \rangle$ for the flat metric on the tangent space of $\mathbb{R}^7\backslash\{0\}$ in $\tilde{p}$. Take $X,\,Y,\,Z\in \mathfrak{X}^1(\mathbb{S}^6)$ and let $\widetilde{X},\,\widetilde{Y}$ and $\widetilde{Z}$ be elements of $\mathfrak{X}(\mathbb{R}^7\backslash\{0\})$ extending the vector fields $X,\,Y$ and $Z$. Then, in the point $N=p$ we have 
    \begin{eqnarray*}
        g_p\big(N_J(X,\,Y),\,JZ\big) & = & -g_p(JN_J(X,\,Y),\,Z)\\
        & = & 4g_p\Big(J\big(\widetilde{J}(\tilde{X}\times \tilde{Y})\big),\,Z\Big)\\
        & = & 4\langle\big(\widetilde{J}\circ\widetilde{J}(\tilde{X}\times \tilde{Y}\big)_p,\,\tilde{Z}_p\rangle\\
        & = & 4\langle -(\tilde{X}\times\tilde{Y})_p^{\perp},\,\tilde{Z}_p\rangle\\
        & = & -4\Big(\langle (\tilde{X}\times\tilde{Y})_p^{\perp},\,\tilde{Z}_p\rangle +\langle (\tilde{X}\times\tilde{Y})_p^{\parallel},\,\tilde{Z}_p\rangle\Big)\\
        & = & -4 \langle (\tilde{X}\times\tilde{Y})_p,\,\tilde{Z}\rangle\\
        & = & -4\widetilde{\omega}_p(\tilde{X}_p,\,\tilde{Y}_p,\,\tilde{Z}_p)\\
        & = & -4\omega_p(X_p,\,Y_p,\,Z_p).
    \end{eqnarray*}
    This shows the claim in $p=N$. By the action of $G_2$ on $\mathbb{S}^6$, we can conclude to the same equality on the whole sphere.
\end{proof}

\begin{thm}\label{automorphisms of omega}
    The connected component of the group of automorphisms of $\mathbb{S}^6$ preserving the $2$-plectic form $\omega$ is $Aut(\mathbb{S}^6,\,\omega)^0= G_2$.
\end{thm}

\begin{proof}
    By Lemma \ref{omega et J} infinitesimal automorphisms preserving $\omega$ have to preserve $J$ as well. Then they preserve $N_J$ by Lemma \ref{NJ G2-invariance}.
    By Proposition \ref{omega=g(NJ...)} these automorphisms have to preserve $g$.
    Thus $Aut(\mathbb{S}^6,\omega)^0\subset Aut(\mathbb{S}^6,g)$ $=O(7,\mathbb{R})$. Since elements of $O(7,\mathbb{R})$ preserving $\omega$ are elements of $G_2$, we conclude that $Aut(\mathbb{S}^6,\omega)^0= G_2$.
\end{proof}

\begin{rem}
     B. Kruglikov shows in \cite{Kruglikov} that the connected component of the group of automorphisms $Aut(\mathbb{S}^6,\,J)^0=G_2$ since he proves that $\mathcal{L}_XJ=0$ implies that $X$ is a fundamental vector field of the $\mathfrak{g}_2$-action. If $J$ is integrable, the Lie algebra of infinitesimal automorphisms of $J$ is $J$-stable but the "strong" non-integrability (compare also Remark \ref{strong non-integrability}) of $J$ underlies the following, seemingly counter-intuitive fact.
\end{rem}

\begin{prop}
    For $\xi\in\mathfrak{g}_2\backslash \{0\}$, the vector field $J(X_{\xi})$ is never a fundamental vector field of the $G_2$-action on $\mathbb{S}^6$.
\end{prop}

\begin{proof}
    Suppose that $X$ and $JX$ are non-zero fundamental vector fields. By Lemma  \ref{omega et J} it follows then that $\mathcal{L}_XJ=0$ and $\mathcal{L}_{JX}J=0$. Thus for all $Y\in\mathfrak{X}^1(\mathbb{S}^6)$,
    \begin{eqnarray*}
        -JN_J(X,Y)=(\mathcal{L}_XJ)(Y)-(\mathcal{L}_{JX}J)(JY)=0,
    \end{eqnarray*}
    i.e. $\forall\,Y\in\mathfrak{X}^1(\mathbb{S}^6),\,N_J(X,\,Y)=0$. With the equality of Proposition \ref{omega=g(NJ...)}, we therefore have
 for all $Y,\,Z\in\mathfrak{X}^1(\mathbb{S}^6)$, $\omega(X,\,Y,\,Z)=-\frac{1}{4}g\big(N_J(X,\,Y)\,,JZ\big)$ $=0$. 
 Since $\omega$ is non-degenerate, we conclude that $X=0$.
\end{proof}

\section{Hamiltonian dynamics on the six-sphere}\label{Hamiltonian dynamics on the six-sphere}

After recalling the notion of (dynamical) Hamilton-de Donder-Weyl equations (compare, e.g. \cite{Wagner&Wurzbacher}), we solve them completely in dimension one and give two different interesting classes of solutions in dimension two (the only possible cases for a $2$-plectic manifold). Notably, we discuss several surprising properties, as the existence of dynamics for a vanishing Hamiltonian function and the fact that the automorphisms of $J$ form a "totally real Lie algebra".

\begin{defn}
    \begin{enumerate}[label=(\roman*)]
        \item Let $(M,\,\omega)$ be an $n$-plectic manifold. A couple $(H,\,X)\in \Omega^{n-k}(M)\times \mathfrak{X}^k(M)$ with $1\leq k\leq n$ is a solution of the \textit{Hamilton-de Donder-Weyl (or HDW) equation} if 
    \begin{eqnarray*}
        \iota_X\omega=(-1)^{n+1-k}dH.
    \end{eqnarray*}
    \item Let $(M,\,\omega)$ be an $n$-plectic manifold, $\Sigma$ a $k$-dimensional manifold with $1\leq k\leq n$ and $H\in \Omega^{n-k}(M)$. A couple $(\gamma,\,\psi)$ 
    with $\gamma\in \mathfrak{X}^k(\Sigma)$ a co-volume and $\psi\::\Sigma\rightarrow M$ a smooth map is called a solution of the \textit{dynamical Hamilton-de Donder-Weyl equation} if for all $x$ in $\Sigma$
    \begin{eqnarray*}
        \iota_{(\psi_*)_x(\gamma_x)}\omega_{\psi(x)}=(-1)^{n+1-k}(dH)_{\psi(x)}.
    \end{eqnarray*}
    \end{enumerate}
    \noindent If $k$ is fixed, we also say that $(H,\,X)$ resp. $(H,\,\gamma,\,\psi)$ is a \textit{solution of the HDW equation} resp. \textit{dynamical HDW equation in 
    dimension $k$ (or equivalently in degree $n$-$k$)}.
\end{defn}

\noindent Thus on $\mathbb{S}^6$ there are only two possibilities\:: $X\in \mathfrak{X}^1(\mathbb{S}^6)$ is a vector field and $H\in \Omega^1(\mathbb{S}^6)$ is a $1$-form, resp. $X\in \mathfrak{X}^2(\mathbb{S}^6)$ is a bivector field and $H\in \Omega^0(\mathbb{S}^6)$ is a function.
\vskip0.2cm

\noindent We first consider the case $(H,\,X)\in\Omega^1(\mathbb{S}^6)\times \mathfrak{X}^1(\mathbb{S}^6)$.

\begin{defn}
    For $\xi\in \mathfrak{g}_2$, let $X_{\xi}$ be the fundamental vector field associated via the $G_2$-action on $\mathbb{S}^6$.
\end{defn}

\begin{prop}
    Let $X\in\mathfrak{X}^1(\mathbb{S}^6)$ and $\omega$ as above. Then $X$ solves the HDW equation if and only if  $X=X_{\xi}$ for a $\xi\in\mathfrak{g}_2$.
\end{prop}

\begin{proof}
    Let $(\alpha,\,X)$ a solution of the HDW equation, i.e. $\iota_X\omega=d\alpha$, then 
    \begin{eqnarray*}
        \mathcal{L}_X\omega & = & d\circ \iota_X\omega +\iota_X(d\omega)\\
        & = & d(d\alpha)=0.
    \end{eqnarray*}
    It follows by Theorem \ref{automorphisms of omega} that $X$ is in Lie\big(Aut($\mathbb{S}^6,\,\omega$)\big)=$\mathfrak{g}_2$.\\
    Let $\theta\in \Omega^2(\mathbb{S}^6)$ be a $G_2$-invariant potential of $\omega$ and $\xi\in\mathfrak{g}_2$. Then
    \begin{eqnarray*}
        \iota_{X_{\xi}}\omega & = & (\iota_{X_{\xi}}\circ d)\theta\\
        & = & \mathcal{L}_{X_{\xi}}\theta-d(\iota_{X_{\xi}}\theta)\\
        & = & d(-\iota_{X_{\xi}}\theta).
    \end{eqnarray*}
    Then $(-\iota_{X_{\xi}}\theta,\,X_{\xi})$ is a solution of the HDW equation.
\end{proof}

\begin{cor} Let $\omega$ be as above on $\mathbb{S}^6$. Then

\begin{enumerate}[label=(\roman*)]
 \item The dimension of the space $\mathfrak{X}_{Ham}(\mathbb{S}^6,\,\omega):=$
       $ \{X\in\mathfrak{X} (\mathbb{S}^6) \, \vert \,  \exists \,  \alpha \in \Omega^1(\mathbb{S}^6)$ $ \text{such that} \,$ $(\alpha,\,X) \text{ is a solution of the HDW equation} \}$
        is finite.
        \item The set of solutions of the HDW equation in dimension one on $(\mathbb{S}^6,\omega)$ is $\{(-\iota_{X_{\xi}}\theta+df,\,X_{\xi})\,|\,\xi\in\mathfrak{g}_2,\,f\in\mathcal{C}^{\infty}(\mathbb{S}^6,\,\mathbb{R})\}$.
    \end{enumerate}
\end{cor}

\begin{prop}
    Given $(H_{\xi},\,X_{\xi})$, a solution of the HDW equation in dimension one, then for all points $p\in\mathbb{S}^6$, the map
    $$\psi^{\xi}_p : \mathbb{R}  \rightarrow  \mathbb{S}^6, \,\, t  \mapsto  e^{t\xi}\cdot p$$
    is a solution of the dynamical HDW equation, i.e. for all $t\in\mathbb{R}$
    \begin{eqnarray*}
        \iota_{(\psi^{\xi}_p)_*\big(\frac{\partial}{\partial t}\big|_t\big)}\omega_{\psi^{\xi}_p(t)}=(dH_{\xi})_{\psi^{\xi}_p(t)}.
    \end{eqnarray*}
    The image of a solution of the dynamical HDW equation is thus the image of an integral curve of a fundamental vector field of the $G_2$-action. If $X_{\xi}(p)=0$, the image is a point, and otherwise the image is a strict immersion of $\mathbb{R}\slash Stab_p(\mathbb{R})$, either diffeomorphic to $\mathbb{S}^1$ or to $\mathbb{R}$. In the latter case, the closure of the orbit is diffeomorphic to $\mathbb{S}^1\times \mathbb{S}^1$.
\end{prop}

\begin{proof}
    All assumptions are clear, up to the characterization of the image of $\psi_p^{\xi}$. Suppose that for $X_{\xi}\in\mathfrak{X}^1(\mathbb{S}^6)$, a solution of the HDW equation, the connected abelian subgroup $\{e^{t\xi}\,|\,t\in\mathbb{R}\}$ of $G_2$ is not $\mathbb{S}^1$, then it is dense in a maximal torus of $G_2$. Such a maximal torus is, of course, diffeomorphic to $\mathbb{S}^1\times \mathbb{S}^1$.
\end{proof}

\begin{rem}
    The leaves of the singular foliation defined by $X_{\xi}$ on $\mathbb{S}^6$ are the images of $\psi^{\xi}_p$ for $p\in\mathbb{S}^6$.
\end{rem}

\noindent Now we consider the case that $(H,\,X)\in\Omega^0(\mathbb{S}^6)\times\mathfrak{X}^2(\mathbb{S}^6)$.

\begin{prop}
    If $\xi,\eta\in\mathfrak{g}_2$ such that $[\xi,\,\eta]=0$ and $H_{(\xi,\eta)}:=-\iota_{X_{\xi}\wedge X_{\eta}}\theta$, then $(H_{(\xi,\eta)},\,X_{\xi}\wedge X_{\eta})$ is a solution of the HDW equation. Furthermore, we have a solution of the dynamical HDW equation
    \begin{eqnarray*}
        \psi_p^{\xi,\eta}(t^1,\,t^2)=\varphi_{t^2}^{X_{\eta}}\circ \varphi_{t^1}^{X_{\xi}}(p)
    \end{eqnarray*}
    for all p in $\mathbb{S}^6$, with $\varphi^Z_t$ denoting the flow of a vector field $Z$ at the time $t$. In this case, the image of $\psi_p^{\xi,\eta}$ is an immersed submanifold diffeomorphic to a torus of dimension $0,\,1$ or $2$ on $\mathbb{S}^6$ or an injectively immersed $\mathbb{R}$ whose closure is a $2$-dimensional torus.
\end{prop}

\begin{proof}
    We have
    \begin{eqnarray*}
        \iota_{X_{\xi}\wedge X_{\eta}}\omega & = & \iota_{X_{\xi}\wedge X_{\eta}}d\theta=\iota_{X_{\eta}}\iota_{X_{\xi}}d\theta\\
        & = & -\iota_{X_{\eta}}d\iota_{X_{\xi}}\theta+\iota_{X_{\eta}}\mathcal{L}_{X_{\xi}}\theta\\
        & = & d\iota_{X_{\eta}}\iota_{X_{\xi}}\theta-\mathcal{L}_{X_{\eta}}\big(\iota_{X_{\xi}}\theta\big)\\
        & = & d\iota_{X_{\eta}}\iota_{X_{\xi}}\theta+\iota_{[X_{\xi},\,X_{\eta}]}\theta\\
        & = & d(\iota_{X_{\eta}}\iota_{X_{\xi}}\theta)=-dH_{(\xi,\eta)}.
    \end{eqnarray*}
    Thus $(-\iota_{X_{\xi}\wedge X_{\eta}}\theta,\,X_{\xi}\wedge X_{\eta})$ is a solution of the HDW equation. Then, for $p\in\mathbb{S}^6$
    \begin{eqnarray*}
        \psi^{\xi,\eta}_p : \mathbb{R}^2 & \rightarrow & \mathbb{S}^6\\
        (t^1,\,t^2) & \mapsto & e^{t^2X_{\eta}}e^{t^1X_{\xi}}\cdot p
    \end{eqnarray*}
    is a solution of the dynamical HDW equation. 
    Since $\xi$ and $\eta$ commute, the connected abelian subgroup $A=\{e^{t^2\eta}e^{t^1\xi}\,|\, t^1,\,t^2\in\mathbb{R}\}$ of $G_2$ is contained in a maximal torus of $G_2$. The rank of $G_2$ being two, this torus is two-dimensional and the closure $\bar{A}$ of $A$ is thus a torus of dimension $0,\,1$ or $2$.\\
    If $A=\bar{A}$, the image of $\psi_p^{\xi,\eta}$ is a submanifold diffeomorphic to $(\mathbb{S}^1)^d$, with $d$ equal or less than the dimension of $A$.\\
    The case $A\subsetneq \bar{A}$ can only rise when the vector space generated by $\xi$ and $\eta$ is one-dimensional and $\bar{A}$ is a maximal torus of $G_2$. The image of $\psi_p^{\xi,\eta}$ is then either $\{p\}$, if $\xi(p)=0=\eta(p)$, or diffeomorphic to $\mathbb{S}^1$, if the stabilizer of $p$ under the action of $A\cong (\mathbb{R},\,+)$ is non-trivial \big(and then isomorphic to $(\mathbb{Z},\,+)$\big), or an immersed submanifold that is diffeomorphic to $\mathbb{R}$ as a manifold and whose closure is $\bar{A}\cdot p$, a two-dimensional product $\mathbb{S}^1\times\mathbb{S}^1$ inside $\mathbb{S}^6$.
\end{proof}

\begin{rem}
    Note that not only solutions of the HDW equations are highly non-unique in dimension two or bigger since a non-vanishing multivector field $Z$ can yield zero upon contraction with $\omega$ but unexpected further non-unicity arrives in the dynamical HDW equations. In the preceding situation, if $\xi\neq 0$ and $\eta=0$ we have $Z=X_{\xi}\wedge X_{\eta}=0$ and $H_{(\xi,\eta)}=0$. Defining $\chi_p(t^1,\,t^2)=p$ we obtain 
    \begin{eqnarray*}
        \big(\chi_p\big)_*\Bigg(\frac{\partial}{\partial t^1}\Bigg|_t\wedge\frac{\partial}{\partial t^2}\Bigg|_t\Bigg)=Z_{\chi_p(t)}=\Big(\psi_p^{\xi,\eta}\Big)_*\Bigg(\frac{\partial}{\partial t^1}\Bigg|_t\wedge\frac{\partial}{\partial t^2}\Bigg|_t\Bigg)
    \end{eqnarray*}
    and then, a fortiori $\chi_p$ and $\psi_p^{\xi,\eta}$ solve the same dynamical HDW equation (with $H_{(\xi,\eta)}$ on the right-hand side). Note also that $\chi_p=\psi_p^{\xi,\eta}$ on the hypersurface $\{t^1=0\}\subset \Sigma=\mathbb{R}^2$.
\end{rem}

\begin{rem}
    With the two previous propositions, we obtain as a corollary the homotopy comoment map of the action of $G_2$ on $\mathbb{S}^6$ of A. M. Miti and L. Ryvkin, compare Lemma 3.8 of \cite{Miti&Ryvkin}.
\end{rem}

\noindent Let us recall that an $n$-plectic form is, in general, only $1$-non-degenerate but not $2$-non-degenerate. Here we can say more. In fact, we even find "dynamics with vanishing Hamiltonian", i.e. solutions $(H,\,Y)\in\Omega^0(\mathbb{S}^6)\times \mathfrak{X}^2(\mathbb{S}^6)$ with $H$ zero or constant.

\begin{prop}
    \begin{enumerate}[label=(\roman*)]
        \item For all $X\in\mathfrak{X}^1(\mathbb{S}^6)$, one has $\iota_{X\wedge JX}\omega=0$, i.e. $(H,\,X\wedge JX)$ for $H$ a constant function is a solution of the HDW equation.
        \item For all $Y\in \mathfrak{X}^2(\mathbb{S}^6)$ decomposable with $\iota_Y\omega=0$, there exists $X\in \mathfrak{X}^1(\{Y\neq 0\})$ such that $Y=\pm X\wedge JX$ on $\{Y\neq 0\}$.
    \end{enumerate}
\end{prop}

\begin{proof}
    \begin{enumerate}[label=(\roman*)]
        \item On $\{X=0\}$ this is obvious. For $p\in\{X\neq 0\}$, we know by Proposition \ref{Lemme technique} that $\ker(\iota_{X_p}\omega)=\big(\big(X_p,\,J(X_p)\big)\big)_{\mathbb{R}}$.
        \item Follows from the last identity of the proof of (i).
    \end{enumerate}
\end{proof}

\begin{prop}\label{commutativity}
    Let $X$ be in $\mathfrak{X}^1(\mathbb{S}^6)$. Then
    \begin{enumerate}[label=(\roman*)]
        \item if $\mathcal{L}_XJ=0$, one has $[X,\,JX]=0$\;;
        \item if $[X,\,JX]=0$ then $D\subset T\mathbb{S}^6$, defined by $D_p=\big(\big(X_p,\,J_p(X_p)\big)\big)_{\mathbb{R}}$ for all p in $\mathbb{S}^6$, is a singular $J$-holomorphic foliation of $\mathbb{S}^6$ whose leaves are images of solutions of the dynamical HDW equations, i.e. for all p in $\mathbb{S}^6$, $\psi_p^{X,\,JX}(t^1,\,t^2)=\varphi_{t^2}^{JX}\circ\varphi_{t^1}^X(p)$.
    \end{enumerate}
\end{prop}

\begin{proof}
    \begin{enumerate}[label=(\roman*)]
        \item We have 
        \begin{eqnarray*}
            [X,\,JX]=\mathcal{L}_X(JX)=(\mathcal{L}_XJ)(X)+J(\mathcal{L}_XX)=0.
        \end{eqnarray*}
        \item The integrable singular distribution $D$ is obviously $J$-holomorphic. Moreover, by the hairy ball theorem, we know that there exist points on $\mathbb{S}^6$ where $X$ vanishes, thus the induced singular foliation has generic leaves of dimension $2$ and singular leaves of dimension $0$.
    \end{enumerate}
\end{proof}

\section*{Acknowledgments}

We thank C. Laurent-Gengoux and L. Ryvkin for discussions related to the content of this article.

\bibliographystyle{spmpsci}
\bibliography{WW-References}

\end{document}